\documentclass[12pt]{amsart}
\usepackage{amssymb}
\usepackage{amsfonts}

\newtheorem{theorem}{Theorem}[section]
\theoremstyle{plain}

\newtheorem{lemma}[theorem]{Lemma}

\newtheorem{proposition}[theorem]{Proposition}

\numberwithin{equation}{section}
\input{tcilatex}

\begin{document}
\title[A Geometric Embedding]{A Geometric Embedding for Standard Analytic
Modules}
\author{Tim Bratten}
\address{Facultad de Ciencias Exactas\\
UNICEN\\
Paraje Arroyo Seco\\
Tandil, Argentina}
\email{bratten@exa.unicen.edu.ar}

\begin{abstract}
In this manuscript we make a general study of the representations realized,
for a reductive Lie group of Harish-Chandra class, on the compactly
supported sheaf cohomology groups of an irreducible finite-rank polarized
homogeneous vector bundle defined in a generalized complex flag space. In
particular, we show that the representations obtained are minimal
globalizations of Harish-Chandra modules and that there exists a whole
number q, depending only on the orbit, such that all cohomologies vanish in
degree less than q. The representation realized on the q-th cohomology group
is called a standard analytic module. Our main result is a geometric proof
that a standard analytic module embeds naturally in an associated standard
module defined on the full flag space of Borel subalgebras. As an
application, we give geometric realizations for irreducible submodules of
some principal series representations in case the group is complex.
\end{abstract}

\maketitle

\section{Introduction}

Even though the irreducible admissible representations of a real reductive
Lie group have been classified for about 30 years, one does not yet have an
explicit realization of these representations. Historically speaking, two of
the first constructions of admissible representations were the parabolically
induced modules and Schmid's realization of the discrete series \cite{thesis}%
, generalizing the Borel-Weil-Bott Theorem \cite{bott}. Developed in the
1980s, the localization technique of Beilinson and Bernstein \cite{algebraic}
gave a general prescription for realizing irreducible Harish-Chandra
modules. However, when the irreducible representation is not obtained as a
standard module, it can be quite difficult to understand the corresponding
geometric realization. Since a specific criteria for the irreducibility of
some standard modules is known \cite{irreducible}, the Beilinson-Bernstein
result does obtain a concrete realization for many irreducible
representations, although the object produced is a Harish-Chandra module and
not a group representation. Using the methods of analytic localization,
Hecht and Taylor gave an analytic version of the Beilinson-Bernstein theory 
\cite{analytic}, including a concrete construction of standard modules that
yields full-fledged group representations globalizing the standard
Beilinson-Bernstein modules. The standard analytic modules defined by Hecht
and Taylor are realized on the compactly supported sheaf cohomology groups
of polarized homogeneous vector bundles. One knows that their construction
is sufficiently general, for example, to provide geometric realizations for
all the tempered representations.

The Hecht-Taylor construction is directly related to the geometry of a full
flag space. In order to realize more irreducible representations, in this
manuscript we consider the natural generalization of their construction to
the setting of an arbitrary flag space. Examples show that, because of
certain geometric conditions on the orbits, the underlying Harish-Chandra
module of this natural generalization is not universally described by the
analog of a standard Beilinson-Bernstein module. In fact, generally
speaking, one knows little about the underlying Harish-Chandra module.
Alternatively, what we propose to establish is an embedding theorem, showing
in general that the standard analytic modules appear naturally as submodules
of certain standard modules defined in a full flag space. In particular, the
standard analytic modules defined in an arbitrary flag space appear as
refinements to certain standard modules occurring in the Beilinson-Bernstein
classification scheme. Our main results along these lines are Theorems 4.1
and 4.3, appearing in Section 4. In particular, the embedding in Theorem 4.1
applies to any finite-rank, irreducible polarized homogeneous vector bundle
defined in a complex flag space. Theorem 4.3 adds a nonvanishing result when
an additional hypothesis is assumed. In the last section, as an example, we
will apply the embedding to realize irreducible submodules of certain
principal series in the case of a complex reductive group.

The paper is organized as follows. The first section is the introduction. In
the second section we introduce the complex flag spaces and the polarized
homogeneous vector bundles. A main result is the construction of the
associated bundle. In the third section we deal with some topological
considerations, introduce the minimal globalization \cite{minimal} and show
that representations obtained on the compactly supported sheaf cohomology
groups of a polarized homogeneous vector bundle are naturally isomorphic to
minimal globalizations of Harish-Chandra modules. We also establish a
vanishing result. In section 4 we prove Theorems 4.1 and 4.3 as mentioned
previously. In section 5 we apply our result to realize irreducible
subrepresentations of some principal series in the case of a complex
reductive group.

\section{Polarized Homogeneous Vector Bundle}

\noindent We let $G_{0}$ denote a reductive Lie group of Harish-Chandra
class with Lie algebra $\mathfrak{g}_{0}$ and complexified Lie algebra $%
\mathfrak{g}$. $G$ denotes the complex adjoint group of $\mathfrak{g}$. This
notation is a bit misleading, since the adjoint representation defines an
isomorphism between the Lie algebra of $G$ and $[\mathfrak{g},\mathfrak{g}]$.

\medskip

\noindent \textbf{Flag Spaces for }$G_{0}$. We define a\emph{\ flag space for%
} $G_{0}$ to be a complex projective, homogeneous $G-$space. The flag spaces
can be constructed as follows. By definition, a\emph{\ Borel subalgebra of }$%
\mathfrak{g}$ is a maximal solvable subalgebra. One knows that $G$ acts
transitively on the set of Borel subalgebras and that the resulting
homogeneous space is a complex projective variety called the\emph{\ full
flag space for} $G_{0}$. A complex subalgebra that contains a Borel
subalgebra is called a\emph{\ parabolic subalgebra of }$\mathfrak{g}$. The
space of $G-$conjugates to a given parabolic subalgebra is a complex
projective variety, and each flag space is realized in this way.

From the construction it follows that points in the flag spaces are
identified with parabolic subalgebras of $\mathfrak{g}$. Let $X$ denote the
full flag space and suppose $Y$ is a flag space. Given $x\in X$ and $y\in Y$
we let $\mathfrak{b}_{x}$ and $\mathfrak{p}_{y}$ denote, respectively, the
corresponding Borel and parabolic subalgebras. For each $x\in X$ it is known
there is a unique $y\in Y$ such that $\mathfrak{b}_{x}\subseteq $ $\mathfrak{%
p}_{y}$. Thus there is a unique $G-$equivariant projection 
\begin{equation*}
\pi :X\rightarrow Y
\end{equation*}%
called the\emph{\ natural projection. }

We note that the quantity of $G_{0}-$orbits in $Y$ is finite \cite{wolf}. In
particular, $G_{0}-$orbits are locally closed in $Y$ and define regular
analytic submanifolds.

\bigskip

\noindent \textbf{The Cartan Dual and Infinitesimal Characters. }Suppose $x$
is a point in the full flag space $X$ and let $\mathfrak{n}_{x}$ denote the
nilradical of the corresponding Borel subalgebra $\mathfrak{b}_{x}$. Put 
\emph{\ } 
\begin{equation*}
\mathfrak{h}_{x}=\mathfrak{b}_{x}/\mathfrak{n}_{x}\;
\end{equation*}%
and let $\mathfrak{h}_{x}^{\ast }$ denote the complex dual of $\mathfrak{h}%
_{x}$. Since the stabilizer of $x$ in $G$ acts trivially on $\mathfrak{h}%
_{x}^{\ast }$ the corresponding $G-$homogeneous holomorphic vector bundle on 
$X$ is trivial. Thus the space $\mathfrak{h}^{\ast }$ of global sections is
isomorphic $\mathfrak{h}_{x}^{\ast }$ via the the evaluation at $x$. We
refer to $\mathfrak{h}^{\ast }$ as the\emph{\ Cartan dual} \emph{for} $%
\mathfrak{g}$. If $\mathfrak{c}$ is a Cartan subalgebra of $\mathfrak{g}$
contained in $\mathfrak{b}_{x}$ then projection of $\mathfrak{c}$ onto $%
\mathfrak{h}_{x}$ coupled with evaluation at $x$ defines an isomorphism of $%
\mathfrak{h}^{\ast }$ onto $\mathfrak{c}^{\ast }$ called\emph{\ }the\emph{\
specialization of }$\mathfrak{h}^{\ast }$ \emph{onto} $\mathfrak{c}^{\ast }$ 
\emph{at} $x$. Using the specializations we can identify a\emph{\ root
subset } 
\begin{equation*}
\Sigma \subseteq \mathfrak{h}^{\ast }
\end{equation*}%
for $\mathfrak{g}$ in $\mathfrak{h}^{\ast }$ and a\emph{\ set of positive
roots} $\Sigma ^{+}\subseteq \Sigma $ where a root is called \emph{positive}
if it is identified with a root of $\mathfrak{c}$ in $\mathfrak{b}_{x}$. The
given subsets are independent of the specialization used to define them. For 
$\alpha \in \Sigma $ and $\lambda \in \mathfrak{h}^{\ast }$ we can also
define the complex values $\overset{\vee }{\alpha }(\lambda )$ and a
corresponding \emph{Weyl group } $W$ \emph{of }$\mathfrak{g}$\emph{\ \ in }$%
\mathfrak{h}^{\ast }$. $\lambda \in \mathfrak{h}^{\ast }$ is called \emph{%
regular }if the stabilizer of $\lambda $ in $W$ is trivial. We say that $%
\lambda $ is \emph{antidominant} provided 
\begin{equation*}
\overset{\vee }{\alpha }(\lambda )\notin \left\{ 1,2,3,\ldots \right\} \text{
\ for each }\alpha \in \Sigma ^{+}.\ 
\end{equation*}

Let $U(\mathfrak{g})$ denote the enveloping algebra of $\mathfrak{g}$ and
let $Z(\mathfrak{g})$ denote the center of $U(\mathfrak{g})$. Then a $%
\mathfrak{g}-$\emph{infinitesimal character} is a homomorphism of algebras 
\begin{equation*}
\Theta :Z(\mathfrak{g})\rightarrow \mathbb{C}\text{.}
\end{equation*}%
A well known result of Harish-Chandra identifies the set of $\mathfrak{g}-$%
infinitesimal characters with the quotient 
\begin{equation*}
\mathfrak{h}^{\ast }/W\text{ .}
\end{equation*}%
Suppose $\Theta $ is an infinitesimal character and $\lambda \in \mathfrak{h}%
^{\ast }$. Then we write $\Theta =W\cdot \lambda $ when Harish-Chandra's
correspondence identifies $\Theta $ with the Weyl group orbit of $\lambda $
in $\mathfrak{h}^{\ast }$.

Let $Y$ be a flag space and suppose $y\in Y$. We let $\mathfrak{p}_{y}$
denote the corresponding parabolic subalgebra and let $\mathfrak{u}_{y}$
denote the nilradical of $\mathfrak{p}_{y}$. The corresponding \emph{\ Levi
quotient} is defined by 
\begin{equation*}
\mathfrak{l}_{y}=\mathfrak{p}_{y}/\mathfrak{u}_{y}\text{ .}
\end{equation*}%
Thus $\mathfrak{l}_{y}$ is a complex reductive Lie algebra. In fact, the
natural projection identifies the Borel subalgebras of $\mathfrak{g}$
contained in $\mathfrak{p}_{y}$ with the Borel subalgebras of $\mathfrak{l}%
_{y}$. In this way, we can identify $\mathfrak{h}^{\ast }$ with the Cartan
dual for $\mathfrak{l}_{y}$. Therefore we have a corresponding \emph{set of
roots} 
\begin{equation*}
\Sigma _{Y}\subseteq \Sigma
\end{equation*}%
for $\mathfrak{h}^{\ast }$ in $\mathfrak{l}_{y}$ and a\emph{\ Weyl group for
the Levi quotient} $W_{Y}\subseteq W$ generated by reflections from the
roots in $\Sigma _{Y}$. As the notation suggests, these subsets are
independent of the point $y$. The\emph{\ positive roots }$\Sigma _{Y}^{+}$
of $\mathfrak{h}^{\ast }$ in $\mathfrak{l}_{y}$ are defined by 
\begin{equation*}
\Sigma _{Y}^{+}=\Sigma _{Y}\cap \Sigma ^{+}.
\end{equation*}%
An element $\lambda \in \mathfrak{h}^{\ast }$ is called \emph{antidominant
for} $Y$ provided there exists an element $w\in W_{Y}$ such that $w\cdot
\lambda $ is antidominant. An equivalent condition is that $\overset{\vee }{%
\alpha }(\lambda )$ not be a positive integer for each $\alpha \in \Sigma
^{+}-\Sigma _{Y}^{+}$. Finally we note that the Harish-Chandra
parametrization identifies set of $\mathfrak{l}_{y}-$infinitesimal
characters with the set of Weyl group orbits 
\begin{equation*}
\mathfrak{h}^{\ast }/W_{Y}\text{ .}
\end{equation*}

\bigskip

\noindent \textbf{Polarized Modules for the Stabilizer. }Suppose $Y$ is a
flag space. Fix $y\in Y$. We have the parabolic subalgebra $\mathfrak{p}_{y}$
with nilradical $\mathfrak{u}_{y}$ and Levi quotient $\mathfrak{l}_{y}$. Let 
$G_{0}[y]$ denote the stabilizer of $y$ in $G_{0}$. Suppose 
\begin{equation*}
\omega :G_{0}[y]\rightarrow GL(V)
\end{equation*}%
is a representation in a finite-dimensional complex vector space $V$. A
compatible representation of $\mathfrak{p}_{y}$ in $V$ is called \emph{a
polarization} if the nilradical $\mathfrak{u}_{y}$ acts trivially. We note
that a polarization need not exist, even if the $G_{0}[y]-$module $V$ is
irreducible. However a polarization is unique when it does exist, because
the compatibility condition assures that two possibly distinct $\mathfrak{l}%
_{y}-$actions agree on a parabolic subalgebra of $\mathfrak{l}_{y}$ and
therefore are identical.

Since a $G_{0}[y]-$invariant subspace of a polarized module need not be
invariant under the corresponding $\mathfrak{l}_{y}-$action, we define \emph{%
a morphism} \emph{of polarized modules} to be a linear map that intertwines
both the $G_{0}[y]$ and $\mathfrak{l}_{y}-$action. Thus the category of
polarized $G_{0}[y]-$modules is nothing but the category of
finite-dimensional $(\mathfrak{l}_{y},G_{0}[y])-$modules.

We briefly consider the structure of an irreducible polarized module $V$.
Since $G_{0}$ is Harish-Chandra class it follows that $V$ has an $\mathfrak{l%
}_{y}-$infinitesimal character. Since there is only one finite-dimensional
irreducible\ $\mathfrak{l}_{y}-$module with a given infinitesimal character,
it follows that $V$ is a direct sum of several copies of an irreducible
highest weight module. In fact, if one makes a linear assumption on the
group $G_{0}$, then an irreducible polarized module is nothing but an
irreducible finite-dimensional $\mathfrak{l}_{y}-$module with compatible $%
G_{0}[y]-$action.

\bigskip

\noindent \textbf{Polarized Homogeneous Vector Bundle. }By definition, a
polarized homogeneous vector bundle is a homogeneous analytic vector bundle
that comes equipped with a canonical way to define the notion of a
restricted holomorphic section. By definition these restricted holomorphic
sections form the corresponding sheaf of polarized sections. For example, in
the case of a trivial line bundle, the polarized sections are exactly the
restricted holomorphic functions.

Let $S$ be a $G_{0}-$orbit and fix $y\in S$. Suppose 
\begin{equation*}
\omega :G_{0}[y]\rightarrow GL(V)
\end{equation*}%
is a polarized module for the stabilizer $G_{0}[y]$. For $\xi \in \mathfrak{p%
}_{y}$ and $v\in V$, we use the notation 
\begin{equation*}
\omega (\xi )v
\end{equation*}%
to indicate the compatible $\mathfrak{p}_{y}-$action. Let 
\begin{equation*}
\begin{array}{c}
\mathbb{V} \\ 
\downarrow \\ 
S%
\end{array}%
\end{equation*}%
be the $G_{0}-$homogeneous vector bundle with fiber $V$ and let 
\begin{equation*}
\phi :G_{0}\rightarrow S\text{ \ \ be the projection }\phi (g)=g\cdot y\text{%
. }
\end{equation*}%
Then a we can identify a section of $\mathbb{V}$ over an open set $%
U\subseteq S$ with a real analytic function 
\begin{equation*}
f:\phi ^{-1}(U)\rightarrow V\text{ \ such that \ }f(gp)=\omega (p^{-1})f(g)%
\text{ \ \ }\forall p\in G_{0}[y].
\end{equation*}%
The section is said to be \emph{polarized} if 
\begin{equation*}
\frac{d}{dt}\mid _{t=0}f(g\exp (t\xi _{1}))+i\frac{d}{dt}\mid _{t=0}f(g\exp
(t\xi _{2}))=-\omega (\xi _{1}+i\xi _{2})f(g)
\end{equation*}%
for all $\xi _{1}$, $\xi _{2}\in \mathfrak{g}_{0}$ such that $\xi _{1}+i\xi
_{2}\in \mathfrak{p}_{y}$.

Let $\mathcal{P}(y,V)$ denote the sheaf of polarized sections and let $%
\mathcal{O}_{Y}\mid _{S}$ be the sheaf of restricted holomorphic functions
on $S$. As a sheaf of $\mathcal{O}_{Y}\mid _{S}-$modules, $\mathcal{P}(y,V)$
is locally isomorphic to $\mathcal{O}_{Y}\mid _{S}\otimes V$ \cite{analytic}%
. The left translation defines a $G_{0}$, and thus a $\mathfrak{g}-$action
on $\mathcal{P}(y,V)$. We will be interested in the representations obtained
on the compactly supported sheaf cohomologies 
\begin{equation*}
H_{\text{c}}^{p}(S,\mathcal{P}(y,V))\text{ \ \ }p=0,1,2,\ldots
\end{equation*}%
One knows that the sheaf cohomologies are dual nuclear Fr\'{e}chet (DNF) $%
\mathfrak{g}-$modules with a compatible $G_{0}-$action, provided certain
naturally defined topologies are Hausdorff \cite{analytic}.

\bigskip

\noindent \textbf{The Associated Bundle. \ }Let $S$ be a $G_{0}-$orbit and
suppose 
\begin{equation*}
\begin{array}{c}
\mathbb{V} \\ 
\downarrow \\ 
S%
\end{array}%
\end{equation*}%
is the homogeneous polarized vector bundle corresponding to an irreducible
polarized module for the stabilizer. We let $\mathcal{P}(\mathbb{V})$ denote
the corresponding sheaf of polarized sections. Let $X$ denote the full flag
space and let 
\begin{equation*}
\pi :X\rightarrow Y
\end{equation*}%
be the natural projection. Put 
\begin{equation*}
N=\pi ^{-1}(S)
\end{equation*}%
the \emph{associated manifold}. In order to study the compactly supported
cohomologies of $\mathcal{P}(\mathbb{V})$, we now define a certain $G_{0}-$%
equivariant bundle over $N$ that comes naturally equipped with a
corresponding polarization.

For $y\in S$ we let $\mathbb{V}_{y}$ denote the fiber of $\mathbb{V}$ over $%
y $. Thus $\mathbb{V}_{y}$ is an irreducible polarized module for the
stabilizer $G_{0}[y]$. Put 
\begin{equation*}
X_{y}=\pi ^{-1}(\left\{ y\right\} )
\end{equation*}%
and let $x\in X_{y}$. Define 
\begin{equation*}
\mathbb{W}_{x}=\frac{\mathbb{V}_{y}}{\mathfrak{n}_{x}\mathbb{V}_{y}}
\end{equation*}%
where $\mathfrak{n}_{x}$ is the nilradical of the Borel subalgebra $%
\mathfrak{b}_{x}$. Thus $\mathbb{W}_{x}$ is an irreducible polarized module
for $G_{0}[x].$ Let $\mathbb{W}$ be the bundle over $N$ given by 
\begin{equation*}
\mathbb{W}=\bigcup\limits_{x\in N}\mathbb{W}_{x}\text{.}
\end{equation*}%
To provide the corresponding trivializations, we define a holomorphic
extension of $\mathbb{W}$ over a neighborhood of $x$ in $X$. Since our
result is local, we may assume that $G_{0}$ is a real form of $G.$ Let $%
B_{x} $ denote the stabilizer of $x$ in $G$. Then the representation of $%
\mathfrak{b}_{x}$ in $\mathbb{W}_{x}$ determines a local holomorphic
representation of $B_{x}$ in $\mathbb{W}_{x}$. The standard construction of
a $G-$homogeneous holomorphic vector bundle limited to an appropriate
neighborhood of the identity in $G$ then provides the desired holomorphic
extension of $\mathbb{W}$ over a neighborhood of $x$ in $X$.

We now consider the polarized sections. For $\xi \in \mathfrak{g}_{0}$ and $%
\sigma $ an analytic section of $\mathbb{W}$ defined on an open set $%
U\subseteq N$ , put 
\begin{equation*}
(\xi \ast \sigma )(x)=\frac{d}{dt}\mid _{t=0}\exp (t\xi )\sigma (\exp (-t\xi
)x),\;\;x\in U
\end{equation*}%
where we use juxtaposition to indicate the action of $G_{0}$ on $\mathbb{W}$%
. This definition determines a unique complex linear action of $\mathfrak{g}$
on the sheaf of analytic sections of $\mathbb{W}$. The section $\sigma $ is
said to be \emph{polarized} if 
\begin{equation*}
(\xi \ast \sigma )(x)=\omega _{x}(\xi )\sigma (x)\text{ \ for each }x\in U%
\text{ and each }\xi \in \mathfrak{b}_{x}
\end{equation*}%
where $\omega _{x}$ indicates the $\mathfrak{b}_{x}-$action on the fiber. In
case 
\begin{equation*}
f:U\rightarrow \mathbb{C}
\end{equation*}%
is a real analytic function, observe that 
\begin{equation*}
(\xi \ast f)(x)=\frac{d}{dt}\mid _{t=0}f(\exp (-t\xi )x),\;\;x\in U\text{ \
and \ }\xi \in \mathfrak{g}_{0}\text{.}
\end{equation*}%
Let $\mathcal{P}(\mathbb{W})$ denote the sheaf of polarized sections and let 
$\mathcal{O}_{X}\mid _{N}$ be the sheaf of restricted holomorphic functions
on $N$. The following lemma will be used to show that $\mathcal{P}(\mathbb{W}%
)$ is a locally free sheaf of $\mathcal{O}_{X}\mid _{N}-$modules.

\begin{lemma}
Suppose $U$ is open in $N$ and let 
\begin{equation*}
f:U\rightarrow \mathbb{C}
\end{equation*}%
be an analytic function. Then $f$ is a restricted holomorphic function if
and only if 
\begin{equation*}
(\xi \ast f)(x)=0\text{ \ \ for each }x\in U\text{ \ and \ }\xi \in 
\mathfrak{b}_{x}\text{.}
\end{equation*}
\end{lemma}

\noindent \textbf{Proof: }If $f$ is a restricted holomorphic function it is
clear that $f$ satisfies the equation in Lemma 2.1. On the other hand,
suppose that $f$ satisfies the given equation. We may assume that $G_{0}$ is
a real form of $G$. To simplify notation, we assume that $f$ is defined on
all of $N$. Let $Q$ be a $G_{0}-$orbit in $N$. We first show that there is a
holomorphic function defined on a neighborhood of $Q$ in $X$ that restricts
to $f$ on $Q$. Fix $x\in Q$ and consider the function $g\mapsto f(gx)$
defined on $G_{0}$. Observe that this function extends to a holomorphic
function $\varphi $ defined on an open set of $G$. For each $g\in G_{0}$ and
each $\xi \in \mathfrak{b}_{gx}$ we have 
\begin{equation*}
(\xi \ast \varphi )(gx)=\frac{d}{dt}\mid _{t=0}\varphi (g\exp (-t\text{Ad}%
(g^{-1})\xi )x)=0.
\end{equation*}%
Thus, for $g$ in an open set of $G$ containing $G_{0}$, we have 
\begin{equation*}
\varphi (gb)=\varphi (g)\text{ \ for each }b\in B_{x}
\end{equation*}%
so that $\varphi $ defines a holomorphic extension of $f$ to a neighborhood
of $Q$ in $X$.

It follows that there exists a holomorphic function on an open set in $X$
which restricts to $f$ on a dense open set of $N$ (the union of all the open 
$G_{0}-$orbits). Choose $y\in \pi (N)$ and put 
\begin{equation*}
X_{y}=\pi ^{-1}(\left\{ y\right\} ).
\end{equation*}%
Define 
\begin{equation*}
\psi :G_{0}\times X_{y}\rightarrow \mathbb{C}\text{ \ by }\psi (g,z)=f(gz)%
\text{. }
\end{equation*}%
We claim that, for each fixed $g$, the function $z\mapsto \varphi (g,z)$ is
holomorphic. By the above considerations, the given function is holomorphic
on a dense open subset of $X_{y}$. But therefore the claim follows, since
the given function is annihilated by the $\overline{\partial }$-operator the
complex manifold $X_{y}$.

Thus $\psi $ extends to a holomorphic function defined on a neighborhood of $%
G_{0}\times X_{y}$ in $G\times X_{y}$. Using the same letter $\psi $ to
denote the holomorphic extension, observe that for each $\xi \in \mathfrak{b}%
_{x}$ and each $(g,z)\in G_{0}\times X_{y}$ we have 
\begin{equation*}
\frac{d}{dt}\mid _{t=0}\psi (g\exp (t\xi ),z)=\frac{d}{dt}\mid _{t=0}\psi
(g,\exp (t\xi )z).
\end{equation*}%
Thus, for each $b\in B_{x}$, we have 
\begin{equation*}
\psi (gb,z)=\psi (g,bz).
\end{equation*}%
It follows that $\psi $ defines a holomorphic extension of $f$ to a
neighborhood of $N$ in $X$. $\blacksquare $

\bigskip

From the previous lemma and the chain rule, it follows that if $\sigma $ is
a local polarized section and $f$ is a restricted polarized section then $%
f\sigma $ is a polarized section. We finish this section with the following
proposition.

\bigskip

\begin{proposition}
Fix $x\in X$. As a sheaf of $\mathcal{O}_{X}\mid _{N}-$modules $\mathcal{P}(%
\mathbb{W})$ is locally isomorphic to $\mathcal{O}_{X}\mid _{N}\otimes 
\mathbb{W}_{x}$.
\end{proposition}

\noindent \textbf{Proof: }Since the result is local, and using the fact $%
\mathbb{W}$ locally extends to a holomorphic vector bundle, we may assume
that $\mathbb{W}$ extends to a $G-$homogeneous vector bundle on all of $X$.
Let 
\begin{equation*}
\Phi :G\rightarrow X\text{ \ be the projection \ }\Phi (g)=gx
\end{equation*}%
and suppose 
\begin{equation*}
\eta :U\rightarrow G
\end{equation*}%
is a corresponding holomorphic section defined on a neighborhood $U$ of $x$
in $X$. We need to check that an analytic section $\sigma $ of $\mathbb{W}$
is polarized if and only if the function 
\begin{equation*}
f(z)=\eta (z)^{-1}\sigma (z)
\end{equation*}%
is holomorphic. Using the differential of the action map 
\begin{equation*}
G\times \mathbb{W}_{z}\rightarrow \mathbb{W}
\end{equation*}%
we define a multiplication by complex scalars in the tangent space to $%
\mathbb{W}$ at point $v\in \mathbb{V}_{z}$. Let $\xi _{1},\xi _{2}\in 
\mathfrak{g}_{0}$ such that $\xi =\xi _{1}+i\xi _{2}\in \mathfrak{b}_{z}$
and put $v=\sigma (z)$. Since 
\begin{equation*}
\frac{d}{dt}\mid _{t=0}\exp (t\xi _{1})v+i\frac{d}{dt}\mid _{t=0}\exp (t\xi
_{2})v=\omega (\xi )v
\end{equation*}%
it follows from the product rule that $\sigma $ satisfies the polarized
condition at $z$ if and only if 
\begin{equation*}
\frac{d}{dt}\mid _{t=0}\sigma (\exp (-t\xi _{1})z)+i\frac{d}{dt}\mid
_{t=0}\sigma (\exp (-t\xi _{2})z)=0.
\end{equation*}%
Since 
\begin{equation*}
\frac{d}{dt}\mid _{t=0}\eta (\exp (-t\xi _{1})z)^{-1}v+i\frac{d}{dt}\mid
_{t=0}\eta (\exp (-t\xi _{2})z)^{-1}v=0
\end{equation*}%
the desired result follows from the product rule as applied to the function $%
f$. $\ \ \blacksquare $

\bigskip

$\mathbb{W}$ is called the \emph{associated bundle on }$N$.

\bigskip

\section{Topological considerations and vanishing.}

In order to show that the compactly supported sheaf cohomologies of a
polarized vector bundle are minimal globalizations, we will apply some of
the results by Hecht and Taylor for certain topological sheaves of $%
\mathfrak{g}-$modules, with compatible $G_{0}-$action \cite{analytic}. We
begin this section by considering some of those results. The following
notations will be adopted. If $S\subseteq Y$ is a locally closed subspace
and $\mathcal{F}$ is a sheaf on $Y$ then $\mathcal{F}\mid _{S}$ denotes the
restriction of $\mathcal{F}$ to $S$. If $\mathcal{F}$ is a sheaf on $S$ then 
$\mathcal{F}^{Y}$ denotes the extension by zero of $\mathcal{F}$ to $Y$.

\bigskip

\noindent \textbf{DNF Sheaves and Analytic Modules. \ }By definition, a 
\emph{DNF space} is a complete, locally convex topological vector space
whose continuous dual, when equipped with the strong topology is a nuclear Fr%
\'{e}chet space. In the obvious fashion, one can define the notion of a DNF
algebra and a corresponding category of DNF modules.

We will be interested in sheaves with DNF structure. A\emph{\ DNF sheaf of
algebras} is a sheaf of algebras such that the space of germs of sections
over each compact subset is a DNF algebra and such that the restriction map
for any nested pair of compact subsets is continuous. When $\mathcal{A}$ is
a DNF sheaf of algebras, then one can define, in the obvious fashion a
corresponding concept of \emph{DNF sheaf of} $\mathcal{A}-$\emph{modules}. A
corresponding \emph{morphism} is a morphism of sheaves of $\mathcal{A}-$%
modules that is continuous for all the given topologies. One can show that
the continuity condition is completely determined by continuity on the
stalks. By the open mapping theorem, a continuous sheaf isomorphism has a
continuous inverse. We note that when $\mathcal{F}$ is a DNF sheaf of $%
\mathcal{A}-$modules on $Y$ and $S\subseteq Y$ is a locally closed subset
then $\left( \mathcal{F}\mid _{S}\right) ^{Y}$ is a DNF sheaf of $\mathcal{A}%
-$modules supported in $S$.

The sheaf $\mathcal{O}$ of holomorphic functions is an important example of
DNF sheaf of algebras. Hecht and Taylor show that when $\mathcal{F}$ is a
sheaf of $\mathcal{O}-$modules whose geometric fibers have countable
dimension and whose stalks are free modules over the stalks of $\mathcal{O}$
then $\mathcal{F}$ carries a unique DNF structure consistent with the
natural DNF structure that exists on stalks. With this structure $\mathcal{F}
$ becomes a DNF sheaf of $\mathcal{O}-$modules.

Suppose $\mathcal{A}$ is a DNF sheaf of algebras on a flag space $Y.$ Since $%
Y$ is compact the space of global sections 
\begin{equation*}
A=\Gamma (Y,\mathcal{A)}
\end{equation*}%
is a DNF algebra. Let $\mathbb{M}_{\text{DNF}}(\mathcal{A})$ denote the
category of DNF sheaves of $\mathcal{A}-$modules and let $\mathbb{M}_{\text{%
DNF}}(A)$ denote the category of DNF $A-$modules. Thus the global sections
defines a functor 
\begin{equation*}
\Gamma :\mathbb{M}_{\text{DNF}}(\mathcal{A})\rightarrow \mathbb{M}_{\text{DNF%
}}(A).
\end{equation*}%
Hecht and Taylor, using a modified \v{C}ech resolution, show that the
category $\mathbb{M}_{\text{DNF}}(\mathcal{A})$ has enough $\Gamma -$acyclic
objects. For $\mathcal{M}\in \mathbb{M}_{\text{DNF}}(\mathcal{A})$ they
prove that the sheaf cohomology groups 
\begin{equation*}
H^{p}(Y,\mathcal{M})\text{ \ }p=1,2,\ldots
\end{equation*}%
will be also be DNF $A-$modules, provided the associated topologies are
Hausdorff. Indeed, their construction is functorial and their work shows
that the topological properties for the sheaf cohomologies are independent
of the choice of $\Gamma -$acyclic resolution, as long as these resolutions
originate from the given category.

We will be interested in the following DNF sheaves of algebras on $Y$ . To
each $\lambda \in \mathfrak{h}^{\ast }$, Beilinson and Bernstein have shown
how to define a twisted sheaf of differential operators (TDO) on the full
flag space $X$ \cite{algebraic}. We let $\mathcal{D}_{\lambda }$ denote the
corresponding TDO with holomorphic coefficients. With our parametrization,
when $\rho $ indicates one-half the sum of the positive roots in $\mathfrak{h%
}^{\ast }$ then $\mathcal{D}_{-\rho }$ is the sheaf of holomorphic
differential operators on $X.$ One knows that $\mathcal{D}_{\lambda }$ is a
DNF sheaf of algebras on $X$. Let 
\begin{equation*}
\pi :X\rightarrow Y
\end{equation*}%
denote the natural projection and let $\pi _{\ast }(\mathcal{D}_{\lambda })$
denote the corresponding direct image in the category of sheaves. Since $X$
and $Y$ are compact, $\pi _{\ast }(\mathcal{D}_{\lambda })$ is a is a DNF
sheaf of algebras on $Y$. If $W_{Y}$ denotes the Weyl group for the Levi
quotient then one knows that $\pi _{\ast }(\mathcal{D}_{w\lambda })\cong \pi
_{\ast }(\mathcal{D}_{\lambda })$ for each $w\in W_{Y}$.

To characterize the sheaf cohomologies of $\ \pi _{\ast }(\mathcal{D}%
_{\lambda })$ we introduce the $\mathfrak{g}-$infinitesimal character $%
\Theta =W\cdot \lambda $. Let $U_{\Theta }$ be the quotient of $U(\mathfrak{g%
})$ by the ideal generated from the kernel of $\Theta $. Using the result of
Beilinson and Bernstein, Hecht and Taylor show that 
\begin{equation*}
\Gamma (Y,\pi _{\ast }(\mathcal{D}_{\lambda }))\cong U_{\Theta }\text{ \ and
\ }H^{p}(Y,\pi _{\ast }(\mathcal{D}_{\lambda }))=0\text{ \ for }p>0.
\end{equation*}

Let $y\in Y$ and suppose $V$ is an irreducible polarized module for the
stabilizer $G_{0}[y]$. Let $\mathcal{P}(y,V)$ denote the corresponding sheaf
of polarized sections defined on the $G_{0}-$orbit $S$ of $y$. Let $%
\mathfrak{l}_{y}$ denote the Levi quotient and let $\lambda \in \mathfrak{h}%
^{\ast }$ be a parameter for the $\mathfrak{l}_{y}-$infinitesimal character
on $V$. Then $\mathcal{P}(y,V)^{Y}$ is a DNF sheaf of $\pi _{\ast }(\mathcal{%
D}_{\lambda })-$modules. Thus the compactly supported sheaf cohomologies 
\begin{equation*}
H_{\text{c}}^{p}(S,\mathcal{P}(y,V))\cong H^{p}(Y,\mathcal{P}(y,V)^{Y})
\end{equation*}%
will be DNF $U_{\Theta }-$modules, provided the corresponding topologies are
Hausdorff.

In order to characterize the $G_{0}-$action on $H_{\text{c}}^{p}(S,\mathcal{P%
}(y,V))$, Hecht and Taylor develop a formalism for DNF sheaves of $\pi
_{\ast }(\mathcal{D}_{\lambda })-$modules with compatible analytic $G_{0}-$%
action. The corresponding objects obtained are referred to as \emph{analytic
sheaves of }$(\pi _{\ast }(\mathcal{D}_{\lambda }),G_{0})-$\emph{modules}.
The sheaf $\mathcal{P}(y,V)^{Y}$ provides the fundamental example of these
sorts of objects. By definition a\emph{\ morphism of analytic sheaves of} $%
(\pi _{\ast }(\mathcal{D}_{\lambda }),G_{0})-$\emph{modules} is a morphism
of DNF $\pi _{\ast }(\mathcal{D}_{\lambda })-$modules which is equivariant
for the $G_{0}-$action.

Suppose $M$ is a DNF $U_{\Theta }-$module equipped with a continuous, linear 
$G_{0}-$action $\omega $. $M$ is called an \emph{analytic } $%
(G_{0},U_{\Theta })-$\emph{module} if:

\medskip

\textbf{(1)} for each $v\in M$ the function 
\begin{equation*}
g\mapsto \omega (g)v
\end{equation*}%
is an analytic function from $G_{0}$ to $M$

\textbf{(2) }the derivative of the $G_{0}-$action agrees with the $\mathfrak{%
g}-$action.

\medskip

\noindent \emph{A morphism of analytic }$(G_{0},U_{\Theta })-$\emph{modules}
is a continuous linear map which is equivariant with respect to the $G_{0}-$%
actions.

It is clear from definitions that the global sections functor determines a
functor from the category of analytic sheaves of $(\pi _{\ast }(\mathcal{D}%
_{\lambda }),G_{0})-$modules to the category of analytic $(G_{0},U_{\Theta
})-$modules. What Hecht and Taylor show with their formalism is that, with
respect to the construction of sheaf cohomology, the analytic group action
is functorial, provided that the sheaf cohomology has a Hausdorff topology.

\bigskip

\noindent \textbf{Minimal Globalization }Let $K_{0}\subseteq G_{0}$ be a
maximal compact subgroup and suppose $M$ is a Harish-Chandra module for $%
(K_{0},\mathfrak{g})$. By definition a globalization $M_{\text{glob}}$ of $M$
is an admissible representation of $G_{0}$ in a complete locally convex
space whose underlying subspace of $K_{0}-$finite vectors is isomorphic to $%
M $ as a $(K_{0},\mathfrak{g})-$module. By now one knows that there exist
several canonical and functorial globalizations. In this article we shall be
interested in realizing the remarkable minimal globalization, $M_{\text{min}%
} $ whose existence was first proved by W. Schmid in \cite{minimal}. When $%
M_{\text{glob}}$ is a globalization of $M$ then the inclusion $M\rightarrow
M_{\text{glob}}$ lifts to a continuous $G_{0}-$equivariant linear inclusion 
\begin{equation*}
M_{\text{min}}\rightarrow M_{\text{glob}}.
\end{equation*}%
This lifting defines an isomorphism of $M_{\text{min}}$ onto the analytic
vectors in a Banach space globalization \cite{minimal2}. In particular, $M_{%
\text{min}}$ is an analytic module. One also knows that the functor of
minimal globalization is exact \cite{minimal2}.

\bigskip

\noindent \textbf{Vanishing Number }It turns out there is a specific
non-negative integer $q$ associated to a $G_{0}-$orbit $S$, such that the
compactly supported sheaf cohomologies of any polarized homogeneous vector
bundle will vanish in degrees less than $q$, but such that the $q-$th
cohomology will not vanish for some polarized homogeneous vector bundles on $%
S$. This value $q$ will be called \emph{the vanishing number} of $S$.

In order to define this number, we introduce some definitions and facts
related to the Matsuki duality \cite{matsuki}. Fix a maximal compact
subgroup $K_{0}\subseteq G_{0}$ and let $K$ be the complexification of $%
K_{0} $. The group $K$ acts algebraically on $Y$ and the choice of $K_{0}$
determines an involutive automorphism 
\begin{equation*}
\theta :\mathfrak{g}\rightarrow \mathfrak{g}
\end{equation*}%
called the \emph{Cartan involution }of $\mathfrak{g}$. A point $y\in Y$ is
called \emph{special} if $\mathfrak{p}_{y}$ contains a Cartan subalgebra $%
\mathfrak{c}$ of $\mathfrak{g}$ such $\theta (\mathfrak{c})=\mathfrak{c}$
and such that that $\mathfrak{g}_{0}\cap \mathfrak{c}$ is a real form for $%
\mathfrak{c}$. A $K-$orbit is said to be \emph{Matsuki dual} to a $G_{0}-$%
orbit if their intersection contains a special point. Matsuki has shown that
each $G_{0}-$orbit (each $K-$orbit) contains a special point and that $K_{0}$
acts transitively on the special points in a given $G_{0}-$orbit (in a given
a $K-$orbit). Thus Matsuki duality establishes a one to one correspondence
between the $G_{0}-$orbits and the $K-$orbits in $Y$. In particular each $%
G_{0}-$orbit has a unique Matsuki dual.

Let $S$ be a $G_{0}-$orbit and let $Q$ be its Matsuki dual. The \emph{%
vanishing number of} $S$ is defined to be the codimension $q$ of $Q$ in $Y$.

\bigskip

\noindent \textbf{The Compactly Supported Cohomologies of a Polarized
Homogeneous Vector Bundle.} We begin this subsection with the following
lemma.

\begin{lemma}
Suppose 
\begin{equation*}
0\rightarrow \mathcal{F}\rightarrow \mathcal{G}\rightarrow \mathcal{H}%
\rightarrow 0
\end{equation*}%
is a short exact sequence of analytic sheaves of $(\pi _{\ast }(\mathcal{D}%
_{\lambda }),G_{0})$-modules on $Y$. Assume that for each $p$, the sheaf
cohomology groups $H^{p}(Y,\mathcal{F})$ and $H^{p}(Y,\mathcal{H})$ are
admissible analytic $G_{0}-$modules, isomorphic to the minimal
globalizations of their underlying Harish-Chandra modules. Then the same is
true of $H^{p}(Y,\mathcal{G})$.
\end{lemma}

\noindent \textbf{\noindent \noindent Proof: }Consider the sequence 
\begin{equation*}
H^{p-1}(Y,\mathcal{H})\rightarrow H^{p}(Y,\mathcal{F})\rightarrow H^{p}(Y,%
\mathcal{G}).
\end{equation*}%
Since a morphism of minimal globalizations has closed range, it follows that
the kernel of the second morphism is closed in $H^{p}(Y,\mathcal{F})$.
Therefore, by applying \cite[Corollary A.11 ]{analytic} to the sequence 
\begin{equation*}
H^{p}(Y,\mathcal{F})\rightarrow H^{p}(Y,\mathcal{G})\rightarrow H^{p}(Y,%
\mathcal{H})
\end{equation*}%
one deduces that $H^{p}(Y,\mathcal{G})$ is Hausdorff. Thus $H^{p}(Y,\mathcal{%
G})$ is an analytic $G_{0}-$module. Hence the desired result follows from 
\cite[Lemma 10.11]{analytic}. $\blacksquare $

\bigskip

With respect to the polarized homogeneous vector bundles in the full flag
space, one has the following result.

\begin{proposition}
Suppose $X$ is the full flag space for $G_{0}$ and $S\subseteq X$ is a $%
G_{0} $-orbit. Choose a point $x\in S$ and suppose $V$ is a polarized module
for $G_{0}[x]$. Let $\mathcal{P}(x,V)$ denote the corresponding sheaf of
polarized sections. Then we have the following. \newline
\textbf{(a)} For each $p=0,1,2,\ldots ,$ the compactly supported sheaf
cohomology group 
\begin{equation*}
H_{\text{c}}^{p}(S,\mathcal{P}(x,V))\text{ }
\end{equation*}%
is an analytic $G_{0}-$module, naturally isomorphic to the minimal
globalization of its underlying Harish-Chandra module.\newline
\textbf{(b)} Let $q$ denote the vanishing number of $S$. Then 
\begin{equation*}
H_{\text{c}}^{p}(S,\mathcal{P}(x,V))=0\text{ \ for }p<q\text{.}
\end{equation*}
\end{proposition}

\noindent \textbf{Proof: \ }When $V$ is a polarized module with
infinitesimal character, then this result is proved in \cite{realizing}.
Since the correspondence 
\begin{equation*}
W\longmapsto \mathcal{P}(x,W)^{X}
\end{equation*}%
defines an exact functor from the category of polarized $G_{0}[x]-$modules
to the category of analytic sheaves of $(\pi _{\ast }(\mathcal{D}_{\lambda
}),G_{0})-$modules, using the previous lemma, the general case follows. $%
\blacksquare $

\medskip

Fix $y\in Y$, let $V$ be an irreducible polarized module for $G_{0}[y]$, put 
$S=G_{0}\cdot y$ and let $N=\pi ^{-1}(S)\subseteq X$ be the associated
manifold. In order to study the compactly supported sheaf cohomologies of
the polarized sections $\mathcal{P}(y,V)$ we introduce the associated vector
bundle $\mathbb{W}$, defined on $N$.

Suppose $y\in S$ is special and let $K[y]$ denote the stabilizer of $y$ in $%
K $. Then the duality between $G_{0}$ and $K-$orbits on $X$ descends to a
1-1 correspondence between $G_{0}[y]$ and $K[y]-$orbits on $X_{y}$ \cite%
{matsuki}. In particular, there exists a unique closed $G_{0}[y]-$orbit on $%
X_{y}$ whose Matsuki dual is the open $K[y]-$orbit. It follows there is a
unique $G_{0}-$orbit $O\subseteq N$ which is closed in $N$. Observe that the
vanishing number of $O$ is equal to the vanishing number of $S$. We call $O$ 
\emph{the associated orbit}.

Applying the previous two results we now establish the following.

\begin{lemma}
Fix $y\in Y$, put $S=G_{0}\cdot y$ and let $V$ be an irreducible polarized
module for $G_{0}[y]$. Let $\mathbb{W}$ be the associated vector bundle
defined on $N=\pi ^{-1}(S)\subseteq X$ and let $\mathcal{P}_{\mathbb{W}}$
denote the corresponding sheaf of polarized sections. Then we have the
following. \newline
\textbf{(a)} For each $p=0,1,2,\ldots ,$ the compactly supported sheaf
cohomology group 
\begin{equation*}
H_{\text{c}}^{p}(N,\mathcal{P}_{\mathbb{W}})\text{ }
\end{equation*}%
is an analytic $G_{0}-$module, naturally isomorphic to the minimal
globalization of its underlying Harish-Chandra module.\newline
\textbf{(b)} Let $q$ denote the vanishing number of $S$ and let $O$ be the
associated orbit in $N$. Then 
\begin{equation*}
H_{\text{c}}^{p}(N,\mathcal{P}_{\mathbb{W}}))=0\text{ \ for }p<q
\end{equation*}%
and there exists a natural, $G_{0}-$equivariant, continuous linear inclusion 
\begin{equation*}
H_{\text{c}}^{q}(N,\mathcal{P}_{\mathbb{W}})\rightarrow H_{\text{c}}^{q}(O,%
\mathcal{P}_{\mathbb{W}}\mid _{O})\text{.}
\end{equation*}
\end{lemma}

\noindent \textbf{Proof:} Put $U=N-O$ and consider the short exact sequence 
\begin{equation*}
0\rightarrow \left( \mathcal{P}_{\mathbb{W}}\mid _{U}\right) ^{X}\rightarrow 
\mathcal{P}_{\mathbb{W}}^{X}\rightarrow \left( \mathcal{P}_{\mathbb{W}}\mid
_{O}\right) ^{X}\rightarrow 0.
\end{equation*}%
Thus (a) follows from the previous two results and the long exact sequence
in cohomology provided we know the result for the sheaf cohomologies of \ $%
\left( \mathcal{P}_{\mathbb{W}}\mid _{U}\right) ^{X}$. By choosing a second $%
G_{0}-$orbit $S_{2}\subseteq U$, that is closed in $U$ and defining $%
U_{2}=U-S_{2}$ we obtain the exact sequence 
\begin{equation*}
0\rightarrow \left( \mathcal{P}_{\mathbb{W}}\mid _{U_{2}}\right)
^{X}\rightarrow \left( \mathcal{P}_{\mathbb{W}}\mid _{U}\right)
^{X}\rightarrow \left( \mathcal{P}_{\mathbb{W}}\mid _{S_{2}}\right)
^{X}\rightarrow 0.
\end{equation*}%
Thus we reduce to the case where $U$ is a $G_{0}-$orbit open in $N$ and the
result follows by Proposition 3.2.

Similar considerations can be used to establish (b). In particular, we only
need to add the observation that if $S^{\shortmid }$ is a $G_{0}-$orbit
contained in $U$ then 
\begin{equation*}
\dim (O)<\dim (S^{\shortmid })
\end{equation*}%
and thus Matsuki duality implies that the vanishing number $q^{\shortmid }$
of $S^{\shortmid }$ satisfies $q^{\shortmid }>q$. Thus, where $U=N-O$, we
obtain 
\begin{equation*}
H^{p}(X,\left( \mathcal{P}_{\mathbb{W}}\mid _{U}\right) ^{X})=0\text{ \ for }%
p<q
\end{equation*}%
and the desired emdedding is obtained from the long exact sequence in
cohomology applied to the first of the previous two short exact sequences. $%
\blacksquare $

\begin{theorem}
Let $Y$ be a flag space for $G_{0}$ and suppose $y\in Y$. Let $V$ be a
polarized $G_{0}[y]$-module and let $\mathcal{P}(y,V)$ denote the
corresponding sheaf of polarized sections. Suppose $S$ is the $G_{0}$-orbit
of $y$ and let $q$ denote the vanishing number of $S$.\ Then we have the
following. \newline
\textbf{(a)} For each $p=0,1,2,\ldots ,$ the compactly supported sheaf
cohomology group 
\begin{equation*}
H_{\text{c}}^{p}(S,\mathcal{P}(y,V))\text{ }
\end{equation*}%
is an admissible analytic $G_{0}$-module, naturally isomorphic to the
minimal globalization of its underlying Harish-Chandra module.\newline
\textbf{(b)} Let $q$ denote the vanishing number of $S$. Then $H_{\text{c}%
}^{p}(S,\mathcal{P}(y,V))$ vanishes for $p<q$.
\end{theorem}

\noindent \textbf{\noindent \noindent Proof: }Using Lemma 3.1 we can reduce
to the case where $V$ is irreducible. Let 
\begin{equation*}
\pi :X\rightarrow Y
\end{equation*}%
be the canonical projection, put $N=\pi ^{-1}(S)$ and let 
\begin{equation*}
\begin{array}{c}
\mathbb{W} \\ 
\downarrow \\ 
N%
\end{array}%
\end{equation*}%
be the associated vector bundle. Let $\mathcal{P}_{\mathbb{W}}$ denote the
corresponding sheaf of polarized sections.

For $x\in X_{y}=\pi ^{-1}(\{y\})$ let $\mathfrak{b}_{x}\subseteq \mathfrak{p}%
_{y}$ be the corresponding Borel subalgebra. Let $\mu $ be the element of $%
\mathfrak{h}^{\ast }$ that specializes to the the $\mathfrak{b}_{x}-$action
on $V/\mathfrak{n}_{x}V$. Thus $\mu $ is the lowest weight in $V$ and the
corresponding $\mathfrak{l}_{y}-$infinitesimal character is given by $%
\lambda =\mu -\rho $. Thus $\mathcal{P}(y,V)^{Y}$ is an analytic sheaf of $%
(\pi _{\ast }(\mathcal{D}_{\lambda }),G_{0})-$modules and $\mathcal{P}_{%
\mathbb{W}}^{X}$ is an analytic sheaf of $(\mathcal{D}_{\lambda },G_{0})-$%
modules.

Observe that the direct image functor $\pi _{\ast }$ induces a functor from
the category of DNF sheaves of $\mathcal{D}_{\lambda }-$modules to the
category of DNF sheaves of $\pi _{\ast }(\mathcal{D}_{\lambda })-$modules.
More specifically, $\pi _{\ast }\mathcal{P}_{\mathbb{W}}^{X}$ is an analytic
sheaf of $(\pi _{\ast }(\mathcal{D}_{\lambda }),G_{0})-$modules. We will
deduce Theorem 3.4 from the following lemma.

\begin{lemma}
Maintain the above notations. Then we have the following. \newline
(a) $\pi _{\ast }\mathcal{P}_{\mathbb{W}}^{X}$ $\cong \mathcal{P}(y,V)^{Y}$
as an analytic sheaf of $(\pi _{\ast }(\mathcal{D}_{\lambda }),G_{0})$%
-modules. \newline
(b) Let $R^{p}\pi _{\ast }$ denote the $p$-th derived functor for the direct
image. Then $R^{p}\pi _{\ast }\left( \mathcal{P}_{\mathbb{W}}^{X}\right) =0$
\ for $p>0$.
\end{lemma}

\noindent \textbf{Proof of the Lemma:} \ Consider the following diagram 
\begin{equation*}
\begin{array}{ccccc}
X_{y} & \overset{i}{\rightarrow } & N & \rightarrow  & X \\ 
\downarrow  &  & \downarrow \pi  &  & \downarrow  \\ 
\left\{ y\right\}  & \rightarrow  & S & \rightarrow  & Y%
\end{array}%
\end{equation*}%
Since the map $\pi $ is proper, by standard sheaf cohomology considerations
it suffices to show that 
\begin{equation*}
\pi _{\ast }\left( \mathcal{P}_{\mathbb{W}}\right) \cong \mathcal{P}(y,V)%
\text{ \ and that }R^{p}\pi _{\ast }\left( \mathcal{P}_{\mathbb{W}}\right) =0%
\text{ \ for }p>0\text{.}
\end{equation*}%
To establish the isomorphism $\pi _{\ast }\left( \mathcal{P}_{\mathbb{W}%
}\right) \cong \mathcal{P}(y,V)$, suppose $z\in S$ and let $\Gamma (X_{z},%
\mathbb{W\mid }_{X_{z}})$ denote the global holomorphic sections of the
holomorphic vector bundle 
\begin{equation*}
\begin{array}{c}
\mathbb{W\mid }_{X_{z}} \\ 
\downarrow  \\ 
X_{z}%
\end{array}%
\end{equation*}%
Let $\mathbb{V}$ denote the polarized homogeneous vector bundle over $S$
corresponding to $V$. Then there is a canonical isomorphism of polarized $%
G_{0}[y]-$modules 
\begin{equation*}
\mathbb{V}_{y}\cong \Gamma (X_{y},\mathbb{W\mid }_{X_{y}})\text{.}
\end{equation*}%
In particular, if $v\in \mathbb{V}_{y}$ and $x\in X_{y}$ then the natural
projection of $v$ into 
\begin{equation*}
\mathbb{W}_{x}=\frac{\mathbb{V}_{y}}{\mathfrak{n}_{x}\mathbb{V}_{y}}
\end{equation*}%
determines a global section $\sigma :X_{y}\rightarrow \mathbb{W}$. Using
this isomorphism, together with the natural $G_{0}-$action on 
\begin{equation*}
\bigcup\limits_{z\in S}\Gamma (X_{z},\mathbb{W\mid }_{X_{z}})
\end{equation*}%
we identify this union with $\mathbb{V}$. Observe that the natural $%
\mathfrak{p}_{z}-$action in $\Gamma (X_{z},\mathbb{W\mid }_{X_{z}})$
corresponds to the $\mathfrak{p}_{z}-$action on $\mathbb{V}_{z}$. Suppose $U$
is an open subset of $S$. We define the isomorphism from $\pi _{\ast }\left( 
\mathcal{P}_{\mathbb{W}}\right) (U)$ to $\mathcal{P}(y,V)(U)$ as follows.
Let $\sigma \in $ $\pi _{\ast }\left( \mathcal{P}_{\mathbb{W}}\right) (U)$
and suppose $z\in U$. Thus the restriction $\sigma \mid _{X_{z}}$ is a
global holomorphic section of the bundle $\mathbb{W\mid }_{X_{z}}$. Put 
\begin{equation*}
\gamma _{\sigma }(z)=\sigma \mid _{X_{z}}\text{ \ for }z\in U\text{.}
\end{equation*}%
Then one checks that $\gamma _{\sigma }$ is a polarized section of $\mathbb{V%
}$ and that the correspondence 
\begin{equation*}
\sigma \mapsto \gamma _{\sigma }
\end{equation*}%
defines the desired isomorphism.

To establish the vanishing, we first note that the sheaf $\mathcal{P}_{%
\mathbb{W}}$ is locally free as a sheaf of $\mathcal{O}_{X}\mid _{N}-$%
modules and that the map $\pi $ is proper. Therefore the desired result can
be deduced from the Borel-Weil-Bott theorem and Grauert's results for the
direct image. In particular, for $z\in S$ let $\mathcal{O}_{X_{z}}$ denote
the sheaf of holomorphic functions on $X_{z}$ and put 
\begin{equation*}
\mathcal{W}(z)=\mathcal{O}_{X_{z}}\tbigotimes\nolimits_{\mathcal{O}_{X}\mid
_{X_{z}}}\mathcal{P}_{\mathbb{W}}\mid _{X_{z}}.
\end{equation*}%
Then $\mathcal{W}(z)$ is naturally isomorphic to the sheaf of holomorphic
sections of the bundle 
\begin{equation*}
\begin{array}{c}
\mathbb{W\mid }_{X_{z}} \\ 
\downarrow \\ 
X_{z}%
\end{array}%
\end{equation*}%
By the Borel-Weil-Bott theorem \cite{bott} 
\begin{equation*}
H^{P}(X_{z},\mathcal{W}(z))=0\text{ \ for }p>0.
\end{equation*}%
Thus the desired result is deduced by standard techniques from Grauert's
Theorem \cite{Grauert}. $\ \blacksquare $

\medskip

To finish the proof of Theorem 3.4, observe that part (a) follows
immediately from Lemma 3.3, the previous lemma and the Leray spectral
sequence. To establish the vanishing from part (b), let $O$ be the
associated orbit in $N$. Then Matsuki duality implies that the vanishing
number of $O$ is equal to $q$. Thus (b) also follows from Lemma 3.3, Lemma
3.5 and the Leray spectral sequence. $\blacksquare $

\medskip

When $V$ is an irreducible polarized $G_{0}[y]-$module and $q$ is the
vanishing number of the $G_{0}-$orbit $S=G_{0}\cdot y$ then the
representation 
\begin{equation*}
H_{\text{c}}^{q}(S,\mathcal{P}(y,V))
\end{equation*}%
will be called \emph{the standard analytic} $G_{0}-$\emph{module induced from%
} $V$.

\section{Geometric Embedding}

Suppose $S\subseteq Y$ is a $G_{0}-$orbit and let $q$ denote the vanishing
number of $S$. Choose $y\in S$ and let $V$ be an irreducible polarized
module for $G_{0}[y]$. We now continue our analysis of the standard analytic
module $H_{\text{c}}^{q}(S,\mathcal{P}(y,V))$.

Let $N=\pi ^{-1}(S)$ be the associated manifold and let $O$ be the
associated orbit in $N$. Thus $O$ has vanishing number $q$. Choose $x\in
O\cap X_{y}$ and let $V/\mathfrak{n}_{x}V$ be the corresponding polarized
module for $G_{0}[x]$. Then we have the \emph{the associated standard
analytic module }$H_{\text{c}}^{q}(O,\mathcal{P}(x,V/\mathfrak{n}_{x}V))$ 
\emph{determined by }$V$.

The geometric embedding of $H_{\text{c}}^{q}(S,\mathcal{P}(y,V))$ in $H_{%
\text{c}}^{q}(O,\mathcal{P}(x,V/\mathfrak{n}_{x}V))$ now follows immediately
from our previous analysis.

\begin{theorem}
There exists a natural, $G_{0}-$equivariant, continuous linear inclusion 
\begin{equation*}
H_{\text{c}}^{q}(S,\mathcal{P}(y,V))\rightarrow H_{\text{c}}^{q}(O,\mathcal{P%
}(x,V/\mathfrak{n}_{x}V))
\end{equation*}
\end{theorem}

\textbf{Proof:} Let 
\begin{equation*}
\begin{array}{c}
\mathbb{W} \\ 
\downarrow \\ 
N%
\end{array}%
\end{equation*}%
be the polarized equivariant vector bundle in $X$ determined from $V$ and
let $\mathcal{P}_{\mathbb{W}}$ denote the sheaf of polarized sections of $%
\mathbb{W}$. Then 
\begin{equation*}
\mathcal{P}_{\mathbb{W}}\mid _{O}\cong \mathcal{P}(x,V/\mathfrak{n}_{x}V)%
\text{.}
\end{equation*}%
Thus the desired result follows from Lemma 3.3, Lemma 3.5 and the Leray
spectral sequence. $\blacksquare $

\bigskip

\noindent \textbf{The Standard Beilinson-Bernstein Modules. }In general, not
much is known about the underlying Harish-Chandra module of $H_{\text{c}%
}^{q}(S,\mathcal{P}(y,V))$. However we can establish non-vanishing result
for $H_{\text{c}}^{q}(S,\mathcal{P}(y,V))$ using information about the
Harish-Chandra module of $H_{\text{c}}^{q}(O,\mathcal{P}(x,V/\mathfrak{n}%
_{x}V))$. In order to state and prove the result, we need to utilize the
Beilinson-Bernstein construction of standard Harish-Chandra modules on the
full flag space, which is known to describe the Harish-Chandra module of a
standard analytic module in that case. We now review some relevant points
surrounding the Beilinson-Bernstein construction. We remark that \emph{many
sheaves in this section are defined on the topological space} $X^{\text{alg}%
} $, where $X^{\text{alg}}$ denotes the full flag space equipped with the $%
G- $invariant Zariski topology.

For $x\in X^{\text{alg}}$ let $K[x]$ denote the stabilizer of $x$ in $K$. By
definition, \emph{a polarized algebraic} $K[x]-$module is a finite
dimensional algebraic $(K[x],\mathfrak{b}_{x})-$module such that the
nilradical $\mathfrak{n}_{x}$ acts trivially. A morphism of polarized
algebraic modules is a linear map, equivariant for both the $K[x]$ and the $%
\mathfrak{b}_{x}-$actions. Suppose that $x$ is special. One knows there
exists a natural equivalence between the categories of polarized algebraic $%
K[x]-$modules and the category of polarized $G_{0}[x]-$modules.

Let $\lambda \in $ $\mathfrak{h}^{\ast }$ and let\ $\mathcal{D}_{\lambda }^{%
\text{alg}}$ be the corresponding twisted sheaf of differential operators
with regular coefficients, defined on the algebraic variety $X^{\text{alg}}$%
. By definition, a \emph{Harish-Chandra sheaf }with parameter $\lambda $ is
a coherent sheaf of $\mathcal{D}_{\lambda }^{\text{alg}}-$modules with
compatible algebraic $K-$action \cite{HMSW}. Let $\mathbb{M}_{\text{coh}}(%
\mathcal{D}_{\lambda }^{\text{alg}},K)$ denote the corresponding category of
Harish-Chandra sheaves. When $\mathcal{M}\in \mathbb{M}_{\text{coh}}(%
\mathcal{D}_{\lambda }^{\text{alg}},K)$ then the sheaf cohomologies 
\begin{equation*}
H^{p}(X^{\text{alg}},\mathcal{M})\text{ \ \ }p=0,1,2,\ldots
\end{equation*}%
are Harish-Chandra modules with infinitesimal character $\Theta =W\cdot
\lambda $.

Let $\mathbb{M}_{\text{HC}}(U_{\Theta },K)$ denote the category of
Harish-Chandra modules for $(\mathfrak{g},K)$ with infinitesimal character $%
\Theta =W\cdot \lambda $. The global sections defines a functor 
\begin{equation*}
\Gamma :\mathbb{M}_{\text{coh}}(\mathcal{D}_{\lambda }^{\text{alg}%
},K)\rightarrow \mathbb{M}_{\text{HC}}(U_{\Theta },K)\text{.}
\end{equation*}%
Beilinson and Bernstein observed that $\Gamma $ has the left adjoint 
\begin{equation*}
\Delta _{\lambda }:\mathbb{M}_{\text{HC}}(U_{\Theta },K)\rightarrow \mathbb{M%
}_{\text{coh}}(\mathcal{D}_{\lambda }^{\text{alg}},K)
\end{equation*}%
given by 
\begin{equation*}
\Delta _{\lambda }(M)=\mathcal{D}_{\lambda }^{\text{alg}}\tbigotimes%
\nolimits_{U_{\Theta }}M
\end{equation*}%
where $M$ is an object in $\mathbb{M}_{\text{HC}}(U_{\Theta },K)$ and $%
U_{\Theta }=\Gamma (X^{\text{alg}},\mathcal{D}_{\lambda }^{\text{alg}})$. In
fact, when $\lambda \in \mathfrak{h}^{\ast }$ is antidominant, Beilinson and
Bernstein have shown that the functor 
\begin{equation*}
\Gamma :\mathbb{M}_{\text{coh}}(\mathcal{D}_{\lambda }^{\text{alg}%
},K)\rightarrow \mathbb{M}_{\text{HC}}(U_{\Theta },K)
\end{equation*}%
is exact and that the adjointness morphism defines an isomorphism between $%
\Gamma \circ \Delta _{\lambda }$ and the identity functor on $\mathbb{M}_{%
\text{HC}}(U_{\Theta },K)$.

Let $\rho $ be one half the sum of the positive roots and let $W$ be a
polarized algebraic $K[x]-$module. When $\mathfrak{b}_{x}$ acts by the
specialization of $\lambda +\rho \in \mathfrak{h}^{\ast }$ then Beilinson
and Bernstein define a corresponding \emph{standard Harish-Chandra sheaf } 
\begin{equation*}
\mathcal{I}(x,W)\in \mathbb{M}_{\text{coh}}(\mathcal{D}_{\lambda }^{\text{alg%
}},K)
\end{equation*}%
\cite{HMSW}. Suppose that $W$ is irreducible. Then there is a unique $%
\lambda \in \mathfrak{h}^{\ast }$such that $\mathfrak{b}_{x}$ acts by $%
\lambda +\rho $. One knows that $\mathcal{I}(x,W)$ has a unique irreducible
Harish-Chandra subsheaf $\mathcal{J}(x,W)\in \mathbb{M}_{\text{coh}}(%
\mathcal{D}_{\lambda }^{\text{alg}},K)$. Put 
\begin{equation*}
I(x,W)=\Gamma (X^{\text{alg}},\mathcal{I}(x,W))\text{ \ and \ }J(x,W)=\Gamma
(X^{\text{alg}},\mathcal{J}(x,W)).
\end{equation*}%
\qquad The pair $(x,W)$ will be called a \emph{classifying data} if $\lambda 
$ is antidominant and if $J(x,W)\neq 0$. When $(x,W)$ is a classifying data%
\emph{\ }then $J(x,W)$ is the unique irreducible $(\mathfrak{g},K)$%
-submodule of $I(x,W)$. An exact criteria for when $(x,W)$ is a classifying
data is known \cite{irreducible}.

\begin{lemma}
Suppose that $(x,W)$ is a classifying data\emph{. }Let $\mathcal{J}(x,W)$
denote the corresponding irreducible Harish-Chandra sheaf in $\mathbb{M}_{%
\text{coh}}(\mathcal{D}_{\lambda },K)$ and let $J(x,W)$ be the global
sections of $\mathcal{J}(x,W)$. Suppose $\mathcal{K}$ is an object from $%
\mathbb{M}_{\text{coh}}(\mathcal{D}_{\lambda },K)$ contained in $\Delta
_{\lambda }(J(x,W))$, whose global sections are zero. Then $\mathcal{K}$ is
contained in the kernel of the adjointness morphism \ 
\begin{equation*}
\Delta _{\lambda }(J(x,W))\rightarrow \mathcal{J}(x,W).
\end{equation*}
\end{lemma}

\noindent \textbf{Proof:} \ Let $\mathcal{N}$ denote the kernel of the
adjointness morphism. Since $J(x,W)\neq 0$, we have the exact sequence 
\begin{equation*}
0\rightarrow \mathcal{N}\rightarrow \Delta _{\lambda }(J(x,W))\rightarrow 
\mathcal{J}(x,W)\rightarrow 0.
\end{equation*}%
Since $\mathcal{J}(x,W)$ is irreducible, $\mathcal{N}$ is a maximal, $K-$%
equivariant subsheaf of coherent $\mathcal{D}_{\lambda }-$submodules in $%
\Delta _{\lambda }(J(x,W))$. Since $\Gamma \circ \Delta _{\lambda }$ is the
identity, the global sections of $\mathcal{N}$ are zero. Consider the sheaf 
\begin{equation*}
\mathcal{M}=\mathcal{K}+\mathcal{N}\subseteq \Delta _{\lambda }(J(x,W)).
\end{equation*}%
Then $\mathcal{M}$ is a $K-$equivariant subsheaf of coherent $\mathcal{D}%
_{\lambda }$-submodules of $\Delta _{\lambda }(J(x,W))$ containing $\mathcal{%
N}$. Since $\mathcal{M}$ has global sections zero, $\mathcal{M}\neq \Delta
_{\lambda }(J(x,W))$. Thus $\mathcal{M}=\mathcal{N}$ and $\mathcal{K}%
\subseteq \mathcal{N}$. $\blacksquare $

\medskip

\noindent \textbf{Nonvanishing for the Geometric Embedding.} Suppose $x\in X$
is special. Then there is a unique $\theta -$stable Cartan subgroup $H_{0}$
of $G_{0}$ contained in $G_{0}[x]$. In particular, if $N_{0}$ indicates the
connected subgroup of $G_{0}[x]$ with Lie algebra $\mathfrak{g}_{0}\cap 
\mathfrak{n}_{x}$ then 
\begin{equation*}
G_{0}[x]=H_{0}\cdot N_{0}\text{ \ with }H_{0}\cap N_{0}=\left\{ 0\right\} .
\end{equation*}%
Thus a finite-dimensional irreducible polarized representation $W$ of $%
G_{0}[x]$ is nothing but a finite-dimensional irreducible representation of $%
H_{0}$ whose derivative we write as $\lambda +\rho \in \mathfrak{h}^{\ast }$%
. The representation $W$ in turn carries a uniquely compatible structure as
an irreducible polarized algebraic $K\left[ x\right] -$module. Let $\mathcal{%
P}(x,W)$ be the corresponding induced sheaf of polarized sections on $%
O=G_{0}\cdot x$ and put 
\begin{equation*}
I(x,W)=\Gamma (X^{\text{alg}},\mathcal{I}(x,W)).
\end{equation*}%
Let $q$ denote the vanishing number of $O$. Then Hecht and Taylor have shown
that $H_{\text{c}}^{q}(O,\mathcal{P}(x,W))$ is naturally isomorphic to the
minimal globalization of $I(x,W)$. Using this result, we are now ready to
establish the nonvanishing clause of the geometric embedding.

\begin{theorem}
Let $S$ be a $G_{0}$-orbit on $Y$ and suppose $q$ denotes the vanishing
number of $S$. Choose a special point $y\in S$ and let $V$ be an irreducible
polarized representation for $G_{0}[y]$. Let $O$ be the associated orbit in $%
X$ and choose $x\in O\cap X_{y}$. Assume $(x,V/\mathfrak{n}_{x}V)$ is a
classifying data. Then $H_{\text{c}}^{q}(S,\mathcal{P}(y,V))$ is non-zero.
In particular, $H_{\text{c}}^{q}(S,\mathcal{P}(y,V))$ is naturally
isomorphic to a topologically closed, non-zero submodule of $H_{\text{c}%
}^{q}(O,\mathcal{P}(x,V/\mathfrak{n}_{x}V))$.
\end{theorem}

\noindent \textbf{Proof:} \ Let $N$ be the associated manifold in $X$ and
put $U=N-O$. Let 
\begin{equation*}
\begin{array}{c}
\mathbb{W} \\ 
\downarrow \\ 
N%
\end{array}%
\end{equation*}%
associated vector bundle and let $\mathcal{P}_{\mathbb{W}}$ denote the
corresponding sheaf of polarized sections. Consider the short exact sequence 
\begin{equation*}
0\rightarrow \left( \mathcal{P}_{\mathbb{W}}\mid _{U}\right) ^{X}\rightarrow
\left( \mathcal{P}_{\mathbb{W}}\right) ^{X}\rightarrow \left( \mathcal{P}_{%
\mathbb{W}}\mid _{O}\right) ^{X}\rightarrow 0.
\end{equation*}%
The proof is by contradiction. Suppose $H_{\text{c}}^{q}(S,\mathcal{P}%
(y,V))=0$. Since 
\begin{equation*}
H_{\text{c}}^{q}(S,\mathcal{P}(y,V))\cong H^{q}(X,\left( \mathcal{P}_{%
\mathbb{W}}\right) ^{X})
\end{equation*}%
we thus obtain a continuous $G_{0}$-equivariant inclusion 
\begin{equation*}
H_{\text{c}}^{q}(O,\mathcal{P}(x,V/\mathfrak{n}_{x}V))\cong H_{\text{c}%
}^{q}(O,\mathcal{P}_{\mathbb{W}}\mid _{O})\rightarrow H_{\text{c}}^{q+1}(U,%
\mathcal{P}_{\mathbb{W}}\mid _{U}).
\end{equation*}%
Put $\Lambda =\left\{ G_{0}-\text{orbits }E\subseteq N:E\text{ has vanishing
number }q+1\right\} $. Then, arguing as in the previous section, we obtain a
continuous $G_{0}-$equivariant inclusion 
\begin{equation*}
H_{\text{c}}^{q+1}(U,\mathcal{P}_{\mathbb{W}}\mid _{U})\rightarrow
\bigoplus_{E\in \Lambda }H_{\text{c}}^{q+1}(E,\mathcal{P}_{\mathbb{W}}\mid
_{E}).
\end{equation*}%
Let $J(x,V/\mathfrak{n}_{x}V)$ denote the unique irreducible Harish-Chandra
submodule of $I(x,V/\mathfrak{n}_{x}V)$ and let $\overline{J(x,V/\mathfrak{n}%
_{x}V)}$ denote the closure of $J(x,V/\mathfrak{n}_{x}V)$ in $H_{\text{c}%
}^{q}(O,\mathcal{P}(x,V/\mathfrak{n}_{x}V))$. By composing the two
inclusions above, one obtains an $E\in \Lambda $ and a continuous $G_{0}-$%
equivariant inclusion 
\begin{equation*}
\overline{J(x,V/\mathfrak{n}_{x}V)}\rightarrow H_{\text{c}}^{q+1}(E,\mathcal{%
P}_{\mathbb{W}}\mid _{E}).
\end{equation*}%
Choose a special point $z\in E\cap X_{y}.$ Then 
\begin{equation*}
H_{\text{c}}^{q+1}(E,\mathcal{P}_{\mathbb{W}}\mid _{E})\cong H_{\text{c}%
}^{q+1}(E,\mathcal{P}(z,V/\mathfrak{n}_{z}V))
\end{equation*}%
so that the underlying Harish-Chandra module of $H_{\text{c}}^{q+1}(E,%
\mathcal{P}_{\mathbb{W}}\mid _{E})$ is the Beilinson-Bernstein module $I(z,V/%
\mathfrak{n}_{z}V).$ We therefore have an inclusion of Harish-Chandra
modules 
\begin{equation*}
J(x,V/\mathfrak{n}_{x}V)\rightarrow I(z,V/\mathfrak{n}_{z}V)\text{. }
\end{equation*}%
Let $\lambda \in \mathfrak{h}^{\ast }$ such that $\mathfrak{b}_{x}$ and $%
\mathfrak{b}_{z}$ act on $V/\mathfrak{n}_{x}V$ and $V/\mathfrak{n}_{z}V$,
respectively by the specialization of $\lambda +\rho $. Put $\Theta =W\cdot
\lambda $ and consider the localization functor 
\begin{equation*}
\Delta _{\lambda }:\mathbb{M}_{\text{HC}}(U_{\Theta },K)\rightarrow \mathbb{M%
}_{\text{coh}}(\mathcal{D}_{\lambda },K).
\end{equation*}%
By our hypothesis, $\lambda $ is antidominant. Let $\mathcal{I}(z,V/%
\mathfrak{n}_{z}V)$ be the standard Harish-Chandra sheaf from $\mathbb{M}_{%
\text{coh}}(\mathcal{D}_{\lambda },K)$ such that 
\begin{equation*}
I(z,V/\mathfrak{n}_{z}V)=\Gamma (X,\mathcal{I}(z,V/\mathfrak{n}_{z}V)).
\end{equation*}%
By the adjointness property, the inclusion of $J(x,V/\mathfrak{n}_{x}V)$ in $%
I(z,V/\mathfrak{n}_{z}V)$ defines a non-zero morphism 
\begin{equation*}
\varphi :\Delta _{\lambda }(J(x,V/\mathfrak{n}_{x}V))\rightarrow \mathcal{I}%
(z,V/\mathfrak{n}_{z}V).
\end{equation*}%
Thus we have the following exact sequence in $\mathbb{M}_{\text{coh}}(%
\mathcal{D}_{\lambda },K)$ 
\begin{equation*}
0\rightarrow \ker (\varphi )\rightarrow \Delta _{\lambda }(J(x,V/\mathfrak{n}%
_{x}V))\overset{\varphi }{\rightarrow }\text{im}(\varphi )\rightarrow 0
\end{equation*}%
where $\ker (\varphi )$ and im$(\varphi )$ denote the kernel and the image
of $\varphi $, respectively. Since the $K-$orbit of $x$ has open
intersection with the fiber $X_{y}$, it follows that $x$ does not belong to
the closure of the $K-$orbit of $z$, thus the stalk im$(\varphi )_{x}$ is
zero. On the other hand, we claim that the stalk of the quotient 
\begin{equation*}
\left( \Delta _{\lambda }(J(x,V/\mathfrak{n}_{x}V))/\ker (\varphi )\right)
_{x}
\end{equation*}%
is not zero. To establish this claim, observe that since $\Gamma \circ
\Delta _{\lambda }$ is the identity, it follows that 
\begin{equation*}
\Gamma (X,\ker (\varphi ))=0.
\end{equation*}%
Let $\mathcal{J}(x,V/\mathfrak{n}_{x}V)$ denote the unique irreducible
object in $\mathbb{M}_{\text{coh}}(\mathcal{D}_{\lambda },K)$ such that 
\begin{equation*}
J(x,V/\mathfrak{n}_{x}V)=\Gamma (X,\mathcal{J}(x,V/\mathfrak{n}_{x}V)).
\end{equation*}%
Therefore, by Lemma 4.2, $\ker (\varphi )$ is contained in the kernel of the
canonical surjection 
\begin{equation*}
\Delta _{\lambda }(J(x,V/\mathfrak{n}_{x}V))\rightarrow \mathcal{J}(x,V/%
\mathfrak{n}_{x}V)\rightarrow 0.
\end{equation*}%
Therefore we have a surjection 
\begin{equation*}
\Delta _{\lambda }(J(x,V/\mathfrak{n}_{x}V))/\ker (\varphi )\rightarrow 
\mathcal{J}(x,V/\mathfrak{n}_{x}V)\rightarrow 0.
\end{equation*}%
This last surjection implies the claim, since $\mathcal{J}(x,V/\mathfrak{n}%
_{x}V)_{x}\neq 0$. $\blacksquare $

\section{Realizing Some Irreducible Subrepresentations of Principal Series
in a Complex Reductive Group.}

In this section we let $G_{0}$ be a connected complex reductive Lie group.
In order to illustrate the embedding of the previous section, we will
realize some irreducible subrepresentations of the principal series.

Let $C_{0}\subseteq G_{0}$ be a Cartan subgroup with Lie algebra $\mathfrak{c%
}_{0}$ and let $\mathfrak{c}\subseteq \mathfrak{g}$ be the corresponding
complexification. Fix a Borel subgroup $B_{0}$ of $G_{0}$ containing $C_{0}$%
. Let $\mathfrak{b}_{0}$ denote the Lie algebra of $B_{0}$ and let $%
\mathfrak{b}\subseteq \mathfrak{g}$ indicate the complexification. Thus the $%
G_{0}-$orbit of $\mathfrak{b}$ in $X$ is closed. Based on this choice there
is a simple 1-1 correspondence between the Weyl group $W_{0}$ of $\mathfrak{c%
}_{0}$ in $\mathfrak{g}_{0}$ and the set of $G_{0}-$orbits in $X$ as
follows. Let $\Sigma (\mathfrak{c})$ indicate the set of roots $\mathfrak{c}$
in $\mathfrak{g}$ and let $\alpha _{0}$ be a root of $\mathfrak{c}_{0}$ in $%
\mathfrak{g}_{0}$. Then the extension of $\alpha _{0}$ to a complex linear
functional $\alpha \in \mathfrak{c}^{\ast }$ is a root of $\mathfrak{c}$ in $%
\mathfrak{g}$. Thus we identify the roots $\Sigma _{0}$ of $\mathfrak{c}_{0}$
in $\mathfrak{g}_{0}$ with a subset of $\Sigma (\mathfrak{c})$. In turn $%
W_{0}$ is identified with a subgroup of the Weyl group of $\mathfrak{c}$ in $%
\mathfrak{g}$. Thus for $w\in W_{0}$ we can define the Borel subalgebra $%
w\cdot \mathfrak{b}$ of $\mathfrak{g}$ containing $\mathfrak{c}$. The
application assigning $w\in W_{0}$ to the $G_{0}-$orbit of $w\cdot \mathfrak{%
b}$ defines the desired 1-1 correspondence. Observe that the dimension of
the $G_{0}-$orbit grows with the length of $w$. In particular, the open $%
G_{0}-$orbit on $X$ corresponds to the longest element of $W_{0}$.

Let 
\begin{equation*}
\chi :H_{0}\rightarrow \mathbb{C}
\end{equation*}%
be a continuous character. Then $\chi $ extends uniquely to a character of
the stabilizer $G_{0}\left[ w\cdot \mathfrak{b}\right] $ of $w\cdot 
\mathfrak{b}$ in $X$ and the derivative $d\chi $ of $\chi $ extends uniquely
to a 1-dimensional complex representation of $w\cdot \mathfrak{b}$. Thus $%
\chi $ determines an irreducible polarized module for $G_{0}\left[ w\cdot 
\mathfrak{b}\right] $. We use the notation 
\begin{equation*}
A(w,\chi )
\end{equation*}%
to indicate the corresponding standard analytic module. When $w$ is the
identity we simply write $A(\chi )$. In this case $A(\chi )$ is the
principal series representation consisting of the space of real analytic
functions 
\begin{equation*}
f:G_{0}\rightarrow \mathbb{C}\text{ \ such that \ }f(gb)=\chi (b^{-1})f(g)%
\text{ \ for }g\in G_{0}\text{ and }b\in B_{0}.
\end{equation*}

Let $d\chi $ indicate the natural extension of the differential of $\chi $
to an element of $\mathfrak{c}^{\ast }$ and identify $d\chi $ with an
element of $\mathfrak{h}^{\ast }$ via the specialization to $\mathfrak{c}%
^{\ast }$ at $w\cdot \mathfrak{b}$. Put 
\begin{equation*}
\lambda =d\chi -\rho
\end{equation*}%
\emph{the shifted differential in }$\mathfrak{h}^{\ast }$. In particular,
the sheaf used to define $A(w,\chi )$ is a sheaf of $\mathcal{D}_{\lambda }-$%
modules.

We will need the following result, which follows directly by an application
of the analytic version of the intertwining functor, established in \cite[%
Proposition 9.10]{analytic} (also consider \cite{irreducible}).

\begin{proposition}
Suppose $\lambda =d\chi -\rho $. Let $\alpha \in \Sigma $ be a simple route
and suppose 
\begin{equation*}
\overset{\vee }{\alpha }(\lambda )\notin \mathbb{Z}\text{.}
\end{equation*}%
Specializing to $\mathfrak{c}^{\ast }$ via the point $w\cdot \mathfrak{b}$
let $\chi _{\alpha }:H_{0}\rightarrow \mathbb{C}$ be the corresponding
character and define 
\begin{equation*}
\gamma :H_{0}\rightarrow \mathbb{C}\text{ \ \ by \ }\chi \cdot \chi _{\alpha
}^{-1}\text{.}
\end{equation*}%
Then there exists a natural isomorphism 
\begin{equation*}
A(w,\chi )\cong A(s_{\alpha }\cdot w,\gamma )
\end{equation*}
\end{proposition}

Observe that, via the identification made above, the shifted differential 
\begin{equation*}
\tau =d\gamma -\rho \in \mathfrak{h}^{\ast }
\end{equation*}%
satisfies the equation $\tau =s_{\alpha }(\lambda )=\lambda -\overset{\vee }{%
\alpha }(\lambda )\alpha $. Hence, when $\lambda $ is antidominant, then so
is $\tau $. In particular, suppose $\beta $ is a root. Then 
\begin{equation*}
\overset{\vee }{\beta }(s_{\alpha }(\lambda ))=2\frac{\left\langle \beta
,s_{\alpha }(\lambda )\right\rangle }{\left\langle \beta ,\beta
\right\rangle }=2\frac{\left\langle s_{\alpha }(\beta ),\lambda
\right\rangle }{\left\langle s_{\alpha }(\beta ),s_{\alpha }(\beta
)\right\rangle }=\overset{\vee }{s_{\alpha }(\beta )}(\lambda ).
\end{equation*}%
Thus the claim is established, since $s_{\alpha }\cdot \Sigma ^{+}=\left(
\Sigma ^{+}-\left\{ \alpha \right\} \right) \cup \left\{ -\alpha \right\} $
and since $\overset{\vee }{\alpha }(\lambda )\notin \mathbb{Z}$.

We assume that $\lambda $ is antidominant and that $A(w,\chi )$ is a
classifying module. Thus $A(w,\chi )$ contains a unique irreducible
submodule $J(w,\chi )$. In this section we treat the case where $w$ is the
identity and consider the possibility of realizing the irreducible submodule 
$J(\chi )$ of $A(\chi )$.

Recall we have identified the roots $\Sigma _{0}$ of $\mathfrak{c}_{0}$ in $%
\mathfrak{g}_{0}$ with subset of $\Sigma (\mathfrak{c})\subseteq \mathfrak{c}%
^{\ast }$. Put 
\begin{equation*}
\Sigma _{0}(\chi )=\left\{ \alpha \in \Sigma _{0}:\overset{\vee }{\alpha }%
(d\chi )\in \mathbb{Z}\right\} .
\end{equation*}%
One knows that $\Sigma _{0}(\chi )$ is a root subspace of $\Sigma _{0}$. The
roots of $\mathfrak{c}$ in $\mathfrak{b}$ define positive systems $\Sigma
_{0}^{+}$ and $\Sigma _{0}^{+}(\chi )=\Sigma _{0}(\chi )\cap \Sigma _{0}^{+}$%
. $\Sigma _{0}(\chi )$ will be called \emph{parabolic }(\emph{at} $\mathfrak{%
b}$) if every simple root in $\Sigma _{0}^{+}(\chi )$ is simple for $\Sigma
_{0}^{+}$.

Suppose $\Sigma _{0}(\chi )$ is parabolic. Put 
\begin{equation*}
\Sigma _{0}(\mathfrak{u})=\Sigma _{0}^{+}-\Sigma _{0}^{+}(\chi )
\end{equation*}%
and define 
\begin{equation*}
S^{+}=\Sigma _{0}^{+}(\chi )\cup -\Sigma _{0}(\mathfrak{u}).
\end{equation*}%
Thus $S^{+}$ is a positive set of roots for $\Sigma _{0}$. In particular,
there exists a unique $w\in W_{0}$ such that 
\begin{equation*}
w\Sigma _{0}^{+}=S^{+}.
\end{equation*}%
Define 
\begin{equation*}
\mathfrak{p}=w\cdot \mathfrak{b}+\dsum\limits_{\alpha \in \Sigma
_{0}^{+}(\chi )}\mathfrak{g}^{-\alpha }
\end{equation*}%
where $\mathfrak{g}^{-\alpha }\subseteq \mathfrak{g}$ is the corresponding
root subspace. Thus $\mathfrak{p}$ is a parabolic subalgebra of $\mathfrak{g}
$ with Levi factor 
\begin{equation*}
\mathfrak{l}=\mathfrak{c}+\dsum\limits_{\alpha \in \Sigma _{0}^{+}(\chi )}%
\mathfrak{g}^{\alpha }+\dsum\limits_{\alpha \in \Sigma _{0}^{+}(\chi )}%
\mathfrak{g}^{-\alpha }\text{ \ \ and nilradical }\mathfrak{u}%
=\dsum\limits_{\alpha \in \Sigma _{0}(\mathfrak{u})}\mathfrak{g}^{-\alpha }.
\end{equation*}

\begin{lemma}
Let $k$ be the number of roots in $\Sigma _{0}(\mathfrak{u})$. Then there is
a sequence of roots 
\begin{equation*}
\alpha _{1},\ldots ,\alpha _{k}\text{ \ with }\alpha _{j}\in \Sigma _{0}(%
\mathfrak{u})
\end{equation*}%
such that $\alpha _{j+1}$ is simple with respect to the Borel subalgebra 
\begin{equation*}
s_{\alpha _{1}}\cdots s_{\alpha _{j}}\mathfrak{b}
\end{equation*}
\end{lemma}

\noindent \textbf{Proof: }In the case where $\Sigma _{0}(\mathfrak{u}%
)=\Sigma _{0}^{+}$ and $\Sigma _{0}^{+}(\chi )=\varnothing $ a reference for
this result can be found in \cite[Lemma 10.6]{analytic}. The more general
case is proved in the same way. $\blacksquare $

\medskip

Observe that 
\begin{equation*}
w=s_{\alpha _{1}}\cdots s_{\alpha _{k}}\text{.}
\end{equation*}%
Define%
\begin{equation*}
\gamma =\chi \cdot \chi _{\alpha _{1}}^{-1}\cdot \cdots \cdot \chi _{\alpha
_{k}}^{-1}.
\end{equation*}%
Then the shifted differential associated to $\gamma $ is antidominant with
respect to the point $w\cdot \mathfrak{b}$. Using the positive system $%
\Sigma _{0}^{+}(\chi )$ with respect to the Levi factor $\mathfrak{l}$, let $%
V$ be the irreducible $\mathfrak{l}-$module with lowest weight $d\gamma $.
Then $V$ is an irreducible polarized representation for the normalizer $G_{0}%
\left[ \mathfrak{p}\right] $ of $\mathfrak{p}$ in $G_{0}.$Let $A(\mathfrak{p}%
,V)$ denote the corresponding standard analytic module. Then we have the
following result.

\begin{theorem}
Assume the principal series representation $A(\chi )$ is a classifying
module and let $J(\chi )\subseteq A(\chi )$ denote the corresponding
irreducible subrepresentation. Suppose $\Sigma _{0}(\chi )$ is parabolic and
define $\mathfrak{p}$ and $V$ as above. Then there is a natural isomorphism 
\begin{equation*}
J(\chi )\cong A(\mathfrak{p},V).
\end{equation*}
\end{theorem}

\noindent \textbf{Proof: }Define $\gamma $ and $w$ as above. Then
Proposition 5.1 and Lemma 5.2 determine a natural isomorphism 
\begin{equation*}
A(\chi )\cong A(w,\gamma ).
\end{equation*}%
Let $N\subseteq X$ be the associated manifold to the $G_{0}-$orbit of $%
\mathfrak{p}$. Observe that $N$ is open in $X$ and that the $G_{0}-$orbit of 
$w\cdot \mathfrak{b}$ is closed in $N$. Therefore we have the embedding%
\begin{equation*}
A(\mathfrak{p},V)\subseteq A(w,\gamma )
\end{equation*}%
given by Theorem 4.1. In addition, $A(\mathfrak{p},V)\neq \left\{ 0\right\} $
by Theorem 4.3. Since the shifted differential for $\gamma $ is
antidominant, and since the $G_{0}-$orbit of $\mathfrak{p}$ is open, it now
follows from the work in \cite{geometry} that $A(\mathfrak{p},V)$ is
irreducible $\blacksquare $

\end{document}